# DISCUSSION OF "BREAKDOWN AND GROUPS" BY P. L. DAVIES AND U. GATHER

By David E. Tyler

*Rutgers, The State University of New Jersey*

The breakdown point has played an important role within robust statistics over the past 25–30 years. A large part of its appeal is that it is easy to explain and easy to understand. It is often interpreted as "the proportion of bad data a statistic can tolerate before becoming arbitrary or meaningless." In this paper Professors Davies and Gather give us a much needed critical look at this seemingly simple concept, and are to be commended for doing so.

Except for some pathological examples, such as a constant functional, one typically presumes that the breakdown point of a statistic or functional cannot be greater than $1/2$. A heuristic justification for this presumption follows from the simple argument "if over half the data is bad, then one cannot distinquish between the good data and the bad data." The authors nicely show that when one more carefully examines this "common sense" argument, then it only appears to be meaningful within a setting with an appropriate group structure. They further challenge the reader to give meaning to this expression outside of such a setting.

One may be able to modify the definition of the breakdown point in some creative way in order to meet the authors' challenge. For example, one could define the breakdown of a statistic to mean it can be made to go to the boundary of all possible values of the statistics. Such a modification would then imply that a constant statistic has breakdown point 0, which is in intuitive agreement with the notion that a constant statistic conveys no information about the data. Such a modification also nullifies the counterexample given by the authors in Section 6.

The intent of this discussion, though, is not to attempt to defend the notion of breakdown outside of the group setting, which may only be possible on a case-by-case basis. Rather, being in general agreement with the arguments made by the authors, the focus of this discussion is to further examine the breakdown point concept within the group setting via two fundamental examples.

**1. Robust principal component vectors.** One approach to robust principal components is to perform a principal component decomposition on a robust covariance matrix rather than on the sample covariance matrix. The asymptotic distribution and influence function of the principal component roots and vectors follow readily from those of the robust covariance matrix; see, for example, Croux and Haesbroeck (2000). However, the breakdown point of the robust covariance matrix has no information regarding the principal component vectors, since the breakdown of a covariance matrix only implies that either the largest root can become arbitrarily large or the smallest root can become arbitrarily small.

To illustrate this point, let $\{\mathbf{x}_1, \ldots, \mathbf{x}_n\}$ represent a sample in $\Re^d$ and let $S_n = \mathcal{Q}\Delta\mathcal{Q}'$ represent the spectral value decomposition of the sample covariance matrix $S_n$. Define $V_n = \mathcal{Q}\Delta^*\mathcal{Q}'$, where $\Delta^* = diagonal\{\lambda_1^2, \ldots, \lambda_d^2\}$ with $\lambda_j$ being a high breakdown point scale statistic for the univariate sample $\{\mathbf{q}_j'\mathbf{x}_1, \ldots, \mathbf{q}_j'\mathbf{x}_n\}$ and where $\mathcal{Q} = [\mathbf{q}_1, \ldots, \mathbf{q}_d]$. That is, we simply replace the eigenvalues of the sample covariance matrix, which correspond to the variances of the sample principal component variables, with robust variances for the sample principal component variables. The resulting statistic $V_n$ has a high breakdown point, namely the breakdown point of the univariate scale statistic used in its definition, whereas the breakdown point of $S_n$ is zero. Both statistics, though, yield the same principal component vectors. So, using a high breakdown point estimate of the covariance matrix for principal components analysis is in itself meaningless, unless one can show some relationship between it and the breakdown of the principal component vectors.

The principal component vector associated with the largest root of the sample covariance matrix can be made arbitrarily close to any given vector by perturbing just one data point. One implicitly assumes this does not occur if a robust covariance matrix is used in place of the sample covariance matrix. Except for contrived examples like the one constructed in the previous paragraph, the proportion of contamination needed to make the largest principal component vector "arbitrary" is likely to be dependent on the separation between the largest root and the other roots of the robust covariance matrix of the uncontaminated data or distribution, as is the case with the influence function. Thus, the best possible bound on the breakdown point is likely dependent on the structure of the uncontaminated data or distribution. As far as this discussant is aware, there are no known results which allow one to quantify this somewhat obvious conjecture, and so it is of interest to see how the results of this paper might apply.

Some meaningful notion of breakdown for a principal component vector is first needed. For the sake of this discussion, consider any robust version of the largest principal component vector, whether or not it is defined via a robust covariance matrix. One usually regards this as an orthogonally equivariant mapping from the data or distribution into the parameter space $\Theta = S^{d-1} = \{\theta \in \Re^d | \theta'\theta = 1\}$, with the parameters $\theta$ and $-\theta$ being viewed as equivalent. Alternatively, $\Theta$ can be taken to be the set of all one-dimensional subspaces. A natural metric on $\Theta$ is the absolute value of the angle between any two elements, that is, $D(\theta_1, \theta_2) = \arccos(|\theta_1'\theta_2|)$, as is used, for example, in Van Aelst and Willems (2004). This metric, however, does not satisfy condition (2.3) of the paper since the largest possible angle is $\pi$. Since $\Theta$ is compact and with no interior points, it is not possible to define a pseudo-metric on $\Theta$ which does satisfy (2.3). An intuitive definition of breakdown, though, can be obtained by simply replacing $\infty$ with $\pi$ in definition (2.4). Breakdown is then naturally interpreted as the proportion of contamination needed for the largest principal component vector to become orthogonal to that obtained from the uncontaminated data.

Since condition (2.3) cannot be made to hold, the results of the paper do not apply here. If one attempts to extend the results of the paper by also replacing $\infty$ with $\pi$ in the definition of G1 given by (3.3), then the set G1 is null. Even if G1 were not null, the crucial step in the proof of Theorem 3.1 is highly dependent on having an unbounded metric. So, it remains an open question as to what types of limits for breakdown are possible for this problem. Using the group equivariance property alone is probably not sufficient for answering this question since the principal component vector associated with any particular root has the same group equivariant property. Some further constraints on what is meant to be a largest principal component vector may be needed to obtain a meaningful bound on the breakdown point. This unsettled question is not specific to principal component vectors, but also applies to any parameter space for which no unbounded metric exists. Such parameter spaces arise naturally, for example, in the areas of directional data analysis and shape theory.

Perhaps some anomaly always arises not only outside of the group setting, but outside of the group setting with unbounded metrics. For principal components, a technicality arises in that it is possible for some data sets or distributions that the largest principal component vector can be any vector within some subspace of dimension $q > 1$, that is, as some $q$-dimensional subspace. For example, any reasonable definition of the largest principal component vector should be any vector at the standard multivariate normal distribution. The complete parameter space then corresponds to the set of all subspaces of $\Re^d$ rather than simply the set of all one-dimensional subspaces. The largest principal component "vector" can still be restricted to equivariant functionals under the group of orthogonal transformations, and the metric $D$ can be extended to this larger parameter space. For example, for bivariate distributions only $\Re^2$ is added to the parameter space and $D$ can be extended by defining $D(\Re^2, \Re^2) = 0$ and $D(\Re^2, \theta) = \pi$. This implies that breakdown occurs when a "well-defined" vector becomes "undefined," that is, becomes $\Re^2$, or vice versa, which is in intuitive agreement with what one thinks of as breakdown. In this setting, though, there exists an orthogonally equivariant functional with breakdown point 1, namely the constant functional $T(P) = \Re^2$, although it is not consistent.

**2. Redescending $M$-estimates of location.** In Section 6 of their paper, the authors note that the meaning of the breakdown point may even be suspect in the well-studied simple univariate location problem. This motivates them to state that "even in the case of equivariance the success of the concept of breakdown point would seem to be more fragile than it is generally supposed." The intent of the discussion here is to elaborate on their remarks by examining in more detail the behavior of the $M$-estimates of location.

For a univariate sample $X^n = \{x_1, \ldots, x_n\}$, an $M$-estimate of location $T(X^n)$ can be defined as a solution to the $M$-estimating equation

$$(2.1) \qquad \sum_{i=1}^{n} \psi\left(\frac{x_i - t}{c}\right) = 0,$$

for some function $\psi$ and tuning constant $c > 0$. A well-known result is that for monotonic, bounded and odd $\psi$-functions, the breakdown point of $T$ is $1/2$. For redescending $\psi$-functions, the breakdown point of the corresponding $M$-estimate is more complicated. Since redescending $\psi$-functions tend to result in multiple solutions to the $M$-estimating equations, it is more convenient to use the alternative definition of an $M$-estimate of location given by

$$(2.2) \qquad T(X^n) = \arg\min_{t} \sum_{i=1}^{n} \rho\left(\frac{x_i - t}{c}\right),$$

for some $\rho$-function. If $\rho$ is differentiable, then the solution to (2.2) also satisfies (2.1) with $\psi = \rho'$. If the function $\rho(r)$ is even, nondecreasing in $|r|$ and bounded, then when $\rho$ is also differentiable the corresponding $\psi$ function is odd and redescends to zero, that is, $\psi(r) \to 0$ as $|r| \to \infty$. Huber (1984) shows the finite-sample contamination breakdown point of such redescending $M$-estimates of location to be

$$(2.3) \qquad \varepsilon^*(T; X^n) = \frac{1 - A(X^n; c)/n}{2 - A(X^n; c)/n},$$

where

$$(2.4) \qquad A(X^n; c) = \min_{t} \sum_{i=1}^{n} \rho\left(\frac{x_i - t}{c}\right),$$

and without loss of generality $\lim_{r \to \infty} \rho(r) = 1$. As the tuning constant $c \to \infty$, the resulting $M$-estimate looks more like the sample mean (provided $\rho$ is differentiable in a neighborhood of zero). Curiously, though, as the tuning constant $c : 0 \to \infty$ one can note that the breakdown point $\varepsilon^*(T; X^n) : 0 \to 1/2$.

A convenient way to gain insight into this formula for the breakdown point (2.3) is by using the relationship between redescending $M$-estimates of location and kernel density estimators. The objective function for an $M$-estimate of univariate location with fixed scale and a kernel density estimate with a given window width, which can be expressed, respectively, as

$$(2.5) \qquad \frac{1}{n} \sum_{i=1}^{n} \rho\left(\frac{x_i - \mu}{c}\right) \quad \text{and} \quad \hat{f}(x) = \frac{1}{nc} \sum_{i=1}^{n} \kappa\left(\frac{x - x_i}{c}\right),$$

have a one-to-one relationship when $\kappa \propto 1 - \rho$. This relationship has been noted, for example, by Chu, Glad, Godtliebsen and Marron (1998). The $M$-estimate of location for a given tuning constant $c$ then corresponds to the mode of the kernel density estimate with window width $c$. The mode of a kernel density estimate based on the Gaussian kernel $\kappa(r) = e^{-r^2/2}/\sqrt{2\pi}$, for example, corresponds to an $M$-estimate of location based on $\rho(r) = 1 - e^{-r^2/2}$, which is referred to as Welsch's $M$-estimate in Splus and MATLAB. Likewise, the mode of a kernel density estimate based on the Epanechnikov kernel corresponds to the skipped-mean.

As an illustrative example, consider the graphs in Figure 1. The data set in the graphs is a simulated data set comprised of 80% standard normal data and 20% normal data with mean 5 and standard deviation 0.1. The three graphs show the kernel density estimates for this data using a Gaussian kernel and the corresponding objective function for the $M$-estimate of location, for increasing values of the tuning constant $c$. In the first graph, the principal mode of the density is centered about the 20% of tightly compacted points. The second graph corresponds to using a larger value of $c$, and the principal mode is located near the mean of the main 80% of the data. The last graph corresponds to using a relatively large value of $c$.

In the first graph, the principal mode would go off to infinity if the more compact 20% of the data were pushed off to infinity, and hence breakdown occurs. Although this $M$-estimate can be made arbitrary under 20% contamination, it is arguable whether the solution is meaningless. On the other hand, in the third graph, although the principal mode is essentially the sample mean, it will not break down even if the 20% were replaced by 45% and allowed to go to infinity since eventually it will fall outside the window width and not impact the central mode. This phenomenon also occurs if the $M$-estimate is made scale equivariant by introducing a robust scale such as the M.A.D.; see Chen and Tyler (2004).

Within the class of redescending $M$-estimates of location, it is arguable whether the breakdown point is more a descriptive property rather than an optimality property. A higher breakdown point redescending $M$-estimate is not necessarily a more desirable estimate. Note also that the nature of breakdown for a redescender differs from that of a monotonic $M$-estimate. The redescender can only break down if a relatively compact cluster of points goes to infinity. If the spread of the "bad" data is greater than that of the "good" data, then a redescending $M$-estimate cannot be broken down even if the "bad" data is in the majority, whereas such contamination would break down a monotonic $M$-estimate.

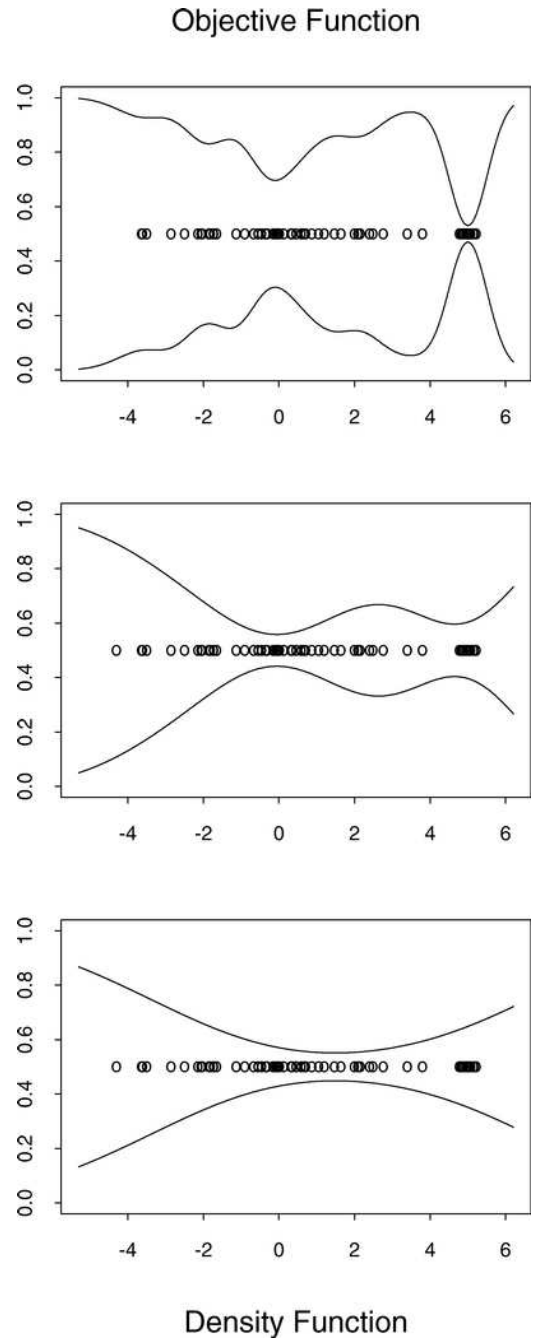

Fig. 1. *An illustration of the relationship between redescending M-estimates and kernel density estimates.*

The above discussion helps illustrate how the simple heuristic interpretation of the breakdown point as "the proportion of bad data a statistical method can tolerate" can be misleading. It has led to some confusion in areas such as computer vision/image understanding. A relatively compact subset of the data may not be considered "bad data" for some applications but rather the data of interest. Instead, a bad data point or "outlier" may be considered a point unlike other data points, and these are the type of bad data points that one may wish to be protected against. At the 2002 ICORS conference in Vancouver, for example, the computer scientist Raymond Ng cleverly noted that a computer scientist's concept of a bad data point can be paraphrased by using the popular Sesame Street phrase "one of these things is not like the others, one of these things does not belong." If one adopts this notion, then the redescending $M$-estimates of location do not break down even under 99% contamination, whereas the monotonic $M$-estimates still have breakdown point $1/2$.

## REFERENCES

Chen, Z. and Tyler, D. E. (2004). On the finite sample breakdown points of redescending $M$-estimates of location. *Statist. Probab. Lett.* **69** 233–242. MR2089000

Chu, C. K., Glad, I. K., Godtliebsen, F. and Marron, J. S. (1998). Edge-preserving smoothers for image processing (with discussion). *J. Amer. Statist. Assoc.* **93** 526–556. MR1631321

Croux, C. and Haesbroeck, G. (2000). Principal component analysis based on robust estimators of the covariance or correlation matrix: Influence functions and efficiencies. *Biometrika* **87** 603–618. MR1789812

Huber, P. J. (1984). Finite sample breakdown points of $M$- and $P$-estimators. *Ann. Statist.* **12** 119–126. MR733503

Van Aelst, S. and Willems, G. (2004). PCA based on multivariate $MM$-estimators with fast and robust bootstrap. Preprint. MR2085872

Department of Statistics
Rutgers, The State University of
  New Jersey
Piscataway, New Jersey 08855
USA
e-mail: dtyler@rci.rutgers.edu